\documentstyle[12pt,amssymb,amstex]{amsart}
\numberwithin{equation}{section}

\def\cb{{\cal B}}

\def\cs{{\cal S}}

\def\ga{{\frak A}}
\def\gb{{\frak B}}

\def\gar{{\frak R}}
\def\gs{{\frak S}}

\def\bc{{\mathbb C}}
\def\bbf{{\mathbb F}}

\def\bm{{\mathbb M}}
\def\bn{{\mathbb N}}

\def\bt{{\mathbb T}}

\def\bz{{\mathbb Z}}

\def\a{\alpha}
\def\b{\beta}
  
\def\d{\delta}

\def\l{\lambda} 

\def\m{\mu}

\def\t{\tau}
\def\f{\varphi} 
  
\def\om{\omega} \def\Om{\Omega}

\newtheorem{thm}{Theorem}[section]

\newtheorem{cor}[thm]{Corollary}
\newtheorem{prop}[thm]{Proposition}
\newtheorem{defin}[thm]{Definition}
\newtheorem{rem}[thm]{Remark}

\def\di{\mathop{\rm d}}

\def\id{\mathop{\rm id}}

\def\idd{{\bf 1}\!\!{\rm I}}

\begin{document}

\title[strict weak mixing]
{strict weak mixing of some $C^*$--dynamical systems
based on free shifts}
\author{Francesco Fidaleo}
\address{Francesco Fidaleo\\
Dipartimento di Matematica\\
II Universit\`{a} di Roma ``Tor Vergata''\\
Via della Ricerca Scientifica, 00133 Roma, Italy}
\email{{\tt fidaleo@@mat.uniroma2.it}}
\author{Farrukh Mukhamedov}
\address{Farrukh Mukhamedov\\
Departamento de Fisica,
Universidade de Aveiro\\
Campus Universitario de Santiago\\
3810-193 Aveiro, Portugal}
\email{{\tt far75m@@yandex.ru},
{\tt farruh@@fis.ua.pt}}

\begin{abstract}
We define a stronger property than unique ergodicity with respect to
the fixed--point subalgebra firstly investigated in \cite{AD}. Such a
property is denoted
as $F$--strict weak mixing ($F$ stands for the Markov projection onto
the fixed--point operator system). Then we show that the free shifts
on the reduced $C^{*}$--algebras of RD--groups,
including the free group on
infinitely many generators, and amalgamated
free product $C^{*}$--algebras,
considered in \cite {AD}, are all strictly weak mixing and not only
uniquely ergodic.
\end{abstract}

\maketitle
\section{Introduction}

Recently, the investigation of the ergodic properties of quantum
dynamical systems had a considerable growth. In quantum setting, the
matter is more complicated than the classical case. For example, some
differences
between classical and quantum situations are pointed out in
\cite{NSZ}. It is then natural to address the study of the possible
generalizations to quantum
case of the various ergodic properties known for classical dynamical
systems.

A very strong ergodic property for a classical system is the unique
ergodicity.
Namely, let $(\Om,T)$ be a classical dynamical systems based on a
compact Haulsdorff space $\Om$ and a homeomorphism $T$ of $\Om$. It
is said to be uniquely ergodic if there exists a unique invariant
Borel measure $\m$ for $T$. It is seen that the ergodic average
${\displaystyle\frac1n\sum_{k=0}^{n-1}f\circ T^{k}}$
converges uniformly to the constant function $\int f\di\m$. The
pivotal 
example of classical uniquely ergodic dynamical system is the
irrational rotations on the unit circle, see e.g. \cite{KFS}. In
quantum setting, the last property is formulated as follows.
Let $(\ga,\a)$ be a $C^{*}$--dynamical systems based on the
$C^{*}$--algebra $\ga$ and the automorphism $\a$. The unique
ergodicity for $(\ga,\a)$ is equivalent (cf. \cite{AD, MT}) to the
norm convergence
\begin{equation}
\label{er}
\lim_{n\to+\infty}\frac1n\sum_{k=0}^{n-1}\a^{n}(a)=E(a)\,,
\end{equation}
where $E$ is the conditional expectation, given by $E=\f(\,\cdot\,)\idd$,
onto the fixed--point
subalgebra of $\a$ consisting of the constant multiples of
the identity. Here, $\f\in\cs(\ga)$ is the unique
invariant state for $\a$. A natural generalization of unique
ergodicity is to require that the ergodic mean in \eqref{er}
converges to a conditional expectation $E$ (necessarily unique)
projecting onto the fixed--point subalgebra $\ga^{\a}$ which, in
general, is
supposed to be nontrivial. This property,
denoted as unique ergodicity with respect to the fixed--point
subalgebra, has been investigated in \cite{AD}. In that paper, it is
proven
that free shifts based on reduced $C^{*}$--algebras of RD--groups
(including the free group on
infinitely many generators), and amalgamated
free product $C^{*}$--algebras, are uniquely ergodic w.r.t. the
fixed--point
subalgebra. This provides nontrivial examples of quantum dynamical
systems based on automorphisms,
exhibiting very strong ergodic properties.

A stronger property than unique
ergodicity, called strict weak
mixing, was investigated in \cite{MT}.
In order to achieve quantum
probability, this was done in the more general situation of
$C^{*}$--dynamical systems $(\ga,T)$, where $T$ is a Markov
(i.e. completely positive and identity--preserving) operator
acting on $\ga$.
The last property simply means that
\begin{equation}
\label{mta}
\lim_{n}\frac1n\sum_{k=0}^{n-1}\big|\psi(T^{n}(a))-\f(a)\big|=0
\end{equation}
for each $\psi\in\cs(\ga)$, $\f$ being the unique invariant state for
$T$. Notice that the
irrational rotations on the unit circle
provides an example of uniquely ergodic dynamical system which is not
strictly weak mixing, see \cite{MT}, Example 2. Other examples
of uniquely ergodic, non strictly weak mixing quantum dynamical systems
can
be easily constructed by
using the algebra of all the $n\times n$ matrices $\bm_{n}(\bc)$, see
\cite{C}.

In the present paper we generalize this mixing--like property to
the situation of \cite{AD}, by considering $C^{*}$--dynamical systems
based on Markov operators. Namely, for the dynamical system
$(\ga,T)$, we require that there exists a linear map
$F:\ga\mapsto\ga$ (necessarily a Markov projection projecting onto the
fixed--point operator system of $T$) such that
\begin{equation}
\label{mta1}
\lim_{n}\frac{1}{n}\sum_{k=0}^{n-1}\big|\psi(T^{k}(x))-\psi(F(x))\big|=0
\end{equation}
for every $\f\in\cs(\ga)$. Such a mixing--like property is denoted
as $F$--strict weak mixing. It is immediate to see that, if there
exists a unique invariant state $\f$ for $T$, then
$F=\f(\,\cdot\,)\idd$ and \eqref{mta1} reduces itself to \eqref{mta}.

Our interest for the above mentioned property is to prove that all the
dynamical systems based on free shifts considered in \cite{AD} are
strictly weak mixing and not merely uniquely ergodic.

\section{terminology, notations and basic results}

Let $\ga$ be a $C^*$-algebra with identity $\idd$. A closed selfadjoint
subspace
$\gar\subset\ga$ containing $\idd$ is said to be an {\it operator system}.
By considering the inclusion $\bm_{n}(\gar)\subset\bm_{n}(\ga)$,
$\bm_{n}(\gar)$ is also an operator system for each $n$. A linear map
$T:\gar\mapsto\gs$ between operator systems is said to be {\it completely
positive} if
$T_{n}:=T\otimes\id_{\bm_{n}}:\bm_{n}(\gar)\mapsto\bm_{n}(\gs)$ is
positive for each $n=1,2,\dots$. It is well--known (cf. \cite{P}) that
$\sup_{n}\|T_{n}\|=T(\idd)$ for completely positive maps.
Let $T:\ga\mapsto\ga$ be completely
positive and identity--preserving, $T$ is called a {\it Markov
operator}.
For such a Markov operator, the fixed--point subspace
$$
\ga^T:=\{x\in\ga\,:\,T(x)=x\}
$$
is an operator system.

Recall that a {\it conditional expectation}
$E:\ga\mapsto\gb\subset\ga$ is a norm--one projection of the
$C^{*}$--algebra $\ga$ onto a $C^{*}$--subalgebra (with the same
identity $\idd$) $\gb$.  The map $E$ is automatically a completely
positive identity--preserving $\gb$--bimodule map, see e.g. \cite{S}.

Let $T$ be a Markov operator. It is seen in \cite{AC} that $\ga^T$ is a
$*$--subalgebra if there exists a faithful invariant state for $T$.
It is readily seen (the same proof as in \cite{AC})
that $\ga^T$ is also a $*$--subalgebra if there exists
a set of invariant states for $T$ which separate the cone of the
positive elements $\ga_{+}$. In general, $\ga^T$ is not a
$*$--subalgebra. At the same way, a Markov projection
$F:\ga\mapsto\ga$ is not
necessarily a conditional expectation, see \cite{CE}, Corollary 7.2.

A (discrete) $C^*$-dynamical system is a pair $\big(\ga,T\big)$
consisting of a $C^*$-algebra and a Markov
operator $T$.

The following theorem is a generalization of Theorem 3.2
of \cite{AD} to our context, see also \cite{MT}, Theorem 3.2 for similar
results.
\begin{thm}
\label{ue}
Let $\big(\ga,T\big)$ be a $C^*$-dynamical system. Then the following
assertions are equivalent.
\begin{itemize}
\item[(i)] Every bounded linear functional on $\ga^T$ has a unique
bounded, $T$--invariant linear extension to $\ga$.
\item[(ii)] Every state on $\ga^T$ has a unique
bounded, $T$--invariant state extension to $\ga$.\footnote{A state
$\f$ on a operator system is a norm one positive linear functional.
In this situation, $\|\f\|=\f(\idd)$.}
\item[(iii)]
$\overline{\ga^{T}+\big\{x-T(x)\,:\,x\in\ga\big\}}=\ga$.
\item[(iv)] The ergodic averages
${\displaystyle\frac{1}{n}\sum_{k=0}^{n-1}T^{k}}$ converge pointwise
in norm.
\item[(v)] The ergodic averages
${\displaystyle\frac{1}{n}\sum_{k=0}^{n-1}T^{k}}$ converge pointwise
in the weak topology to a linear map $F:\ga\mapsto\ga$.\footnote{It is
seen in the proof that such a linear map $F$ is
necessarily a Markov projection onto $\ga^{T}$, satisfying $TF=F=FT$.}
\item[(vi)]
$\ga^{T}+\overline{\big\{x-T(x)\,:\,x\in\ga\big\}}=\ga$.
\end{itemize}

Furthermore, if one (and hence all) of the above statements holds, then
there
exists a unique Markov projection $F$ of $\ga$ onto $\ga^{T}$. It is
given by
\begin{equation}
\label{mp}
F(x)=\lim_{n}\frac{1}{n}\sum_{k=0}^{n-1}T^{k}(x)\,.
\end{equation}
\end{thm}
\begin{pf} (iv)$\Rightarrow$(v): If the ergodic averages
${\displaystyle\frac{1}{n}\sum_{k=0}^{n-1}T^{k}}$ converge in
norm, they define a linear map on $\ga$ given in \eqref{mp}. Then (v)
easily follows. We then prove (v)$\Rightarrow$(ii)
as the remaining implications follow the same lines of
Theorem 3.2 of \cite{AD}.

It is readily seen
that if there exists a linear map $F:\ga\mapsto\ga$ such that, for
$x\in\ga$,
$$
\mathop{\rm w\!-\!lim}_{n}\frac{1}{n}\sum_{k=0}^{n-1}T^{k}(x)=F(x)\,,
$$
then $F$ is automatically a completely positive and
identity--preserving, hence bounded.
In addition, if $f\in\ga^{*}$,
\begin{align*}
&f(F(T(x)))=\lim_{n}\frac{1}{n}\sum_{k=0}^{n-1}f(T^{k+1}(x))\\
=&f(F(x))=\lim_{n}\frac{1}{n}\sum_{k=0}^{n-1}f\circ T(T^{k}(x))
=f(T(F(x)))\,.
\end{align*}
This leads to $TF=F=FT$. Thus,
$$
f(F^{2}(x))=\lim_{n}\frac{1}{n}\sum_{k=0}^{n-1}f(T^{k}(F(x)))
=\lim_{n}\frac{1}{n}\sum_{k=0}^{n-1}f(F(x))
\equiv f(F(x))\,,
$$
that is $F^{2}=F$. Now, if $x=F(x)$ then $T(x)=T(F(x))=F(x)=x$. If
$x=T(x)$ then
$$
f(F(x))=\lim_{n}\frac{1}{n}\sum_{k=0}^{n-1}f(T^{k}(x))=f(x)\,,
$$
that is $x=F(x)$. Namely, $F$ projects onto $\ga^{T}$.

Let now $\f_{j}$,
$j=1,2$ invariant state extensions of the state $\om$ on $\ga^{T}$.
Then
$$
\f_{j}(x)=\lim_{n}\frac{1}{n}\sum_{k=0}^{n-1}\f_{j}(T^{k}(x))
=\f_{j}(F(x))=\om(F(x))\,,
$$
that is any state on $\ga^{T}$ has a unique invariant state extension
on $\ga$.
\end{pf}
\begin{defin}
The	$C^*$-dynamical system $\big(\ga,T\big)$ is said to be
$F$--uniquely ergodic if one of the equivalent properties
$\text(i)$--$\text(vi)$ of Theorem \ref{ue} holds true.\footnote{In
\cite{AD},
the analogous property relative to $C^*$-dynamical systems based on
automorphisms is denoted as {\it unique ergodicity w.r.t. its
fixed--point subalgebra}.}

The	$C^*$-dynamical system $\big(\ga,T\big)$ is said to be
$F$--strictly weak mixing if there exists a linear map
$F:\ga\mapsto\ga$ such that
\begin{equation}
\label{mp1}
\lim_{n}\frac{1}{n}\sum_{k=0}^{n-1}\big|\f(T^{k}(x))-\f(F(x))\big|=0
\end{equation}
whenever $\f\in\cs(\ga)$.
\end{defin}
\begin{prop}
\label{uesm}
If the	$C^*$-dynamical system $\big(\ga,T\big)$ is $F$--strictly weak
mixing, then it is $F$--uniquely ergodic.
\end{prop}
\begin{pf}
Let $\f\in\cs(\ga)$. We get
$$
\bigg|\frac{1}{n}\sum_{k=0}^{n-1}\big(\f(T^{k}(x))-\f(F(x))\big)\bigg|
\leq\frac{1}{n}\sum_{k=0}^{n-1}\big|\f(T^{k}(x))-\f(F(x))\big|\to0
$$
whenever $n\to+\infty$, as $\big(\ga,T\big)$ is $F$--strictly weak
mixing. By using the Jordan decomposition of bounded linear
functionals (cf. \cite{T}), we
conclude that (v) of Theorem \ref{ue} is satisfied.
\end{pf}
\begin{rem}
\label{rem}
By taking into account Proposition \ref{uesm} and Theorem \ref{ue}, the
map $F$ in
\eqref{mp1} is a
Markov projection projecting onto $\ga^{T}$.
\end{rem}

In many interesting situations, the ergodic behavior of dynamical
systems are connected with some spectral properties, see e.g.
\cite{C, L, NSZ}. It is not possible to extend such results to the 
full generality. However, a $F$--strictly weak mixing map $T$ cannot have
eigenvalues on the unit circle $\bt$ except $z=1$.
Namely, for $z$ in $\bc$ denote
$$
\ga_z=\{x\in\ga: T(x)=zx\}.
$$
Of course, $\ga_1=\ga^{T}$. Furthermore,
\begin{prop}
\label{mix2}
Let $\big(\ga,T\big)$ be a $F$--strictly weak mixing
$C^*$-dynamical system. Then $z\in\bt\backslash\{1\}$ implies
$\ga_z=\{0\}$.
\end{prop}
\begin{pf}
Assume that $T(x_0)=zx_0$ for some $z\neq 1$.
Then $F(x_0)=F(T(x_0))=zF(x_0)$ which means $F(x_0)=0$. In addition,
from the $F$-strict weak
mixing we infer
\begin{align*}
0=&\lim_{n}\frac{1}{n}\sum_{k=0}^{n-1}\big|\f(T^k(x_0))-\f(F(x_0))\big|
=\lim_{n}\frac{1}{n}\sum_{k=0}^{n-1}\big|z^{k}\f(x_0)\big|\\
=&\lim_{n}\frac{1}{n}\sum_{k=0}^{n-1}\big|\f(x_0)\big|
= |\f(x_0)|\,.
\end{align*}
Namely, $\f(x_0)=0$ for every $\f\in\cs(\ga)$, hence $x_0=0$.
\end{pf}

Finally, we recall the following result relative to the bounded
sequences which are weakly mixing to zero. For the definitions of
(lower) density of a subset of the natural numbers, or relatively
dense sequences of natural numbers, we refer the reader to \cite{Z}.
\begin{thm} (cf. \cite{Z}, Theorem 2.3)
\label{z}
Let $\{x_{n}\}_{n\geq1}$ be a bounded sequence in the Banach space $X$.
The following assertions are equivalent.
\begin{itemize}
\item[(i)]
${\displaystyle\lim_{n}\,\,\,\sup\bigg\{\frac{1}{n}\sum_{k=1}^{n}\big|f(x_k)\big|\,:\,
f\in X^{*}\,,\,\|f\|\leq1\bigg\}=0}$.
\item[(ii)]
${\displaystyle\lim_{n}\bigg\|\frac{1}{n}\sum_{j=1}^{n}x_{k_{j}}\bigg\|=0}$
for each sequence $k_{1}<k_{2}<\cdots$ of strictly positive lower
density.
\item[(iii)]
${\displaystyle\lim_{n}\bigg\|\frac{1}{n}\sum_{j=1}^{n}x_{k_{j}}\bigg\|=0}$
for each relatively dense sequence\\
$k_{1}<k_{2}<\cdots$.
\end{itemize}
\end{thm}

\section{strict weak mixing of lenght--preserving automorphisms
of RD--groups}

In \cite{AD}, it has been proved that some automorphisms of the
reduced $C^{*}$--algebra of RD--groups are $E$--uniquely
ergodic.\footnote{The
RD--groups are defined and studied in \cite{J}. Notice that
the RD--groups include the Gromov hyperbolic groups, see \cite{dH}.}
Here, we prove that such automorphisms are $E$--strictly weak mixing.
\begin{prop}
\label{ue1}
Let $\b$ be a lenght--preserving automorphism of a RD--group $G$ for
the lenght--function $L$, such
that its orbits are infinite or singletons. Then the automorphism $\a$
induced by $\b$ on
$C^{*}_{r}(G)$ is $E$--strictly weak mixing.
\end{prop}
\begin{pf}
Let $H:=\{g\in G\,:\,\b(g)=g\}$.
As $\a$ is $E$--uniquely ergodic (cf. \cite{AD}, Proposition 3.5),
the pointwise limit in norm
$$
E:=\frac{1}{n}\sum_{k=1}^{n}\a^{k}
$$
exists and gives rise a conditional expectation projecting onto
the fixed--point algebra $C^{*}_{r}(H)\subset C^{*}_{r}(G)$.
By a standard density argument, it is enough to prove that the
sequence $\{\a^{n}(\l_{g})\}_{n\geq1}$ is weakly mixing to zero
whenever $\b(g)\neq g$, that is
$$
\frac{1}{n}\sum_{k=1}^{n}\big|f(\a^{k}(\l_{g}))\big|\to0
$$
for each $f\in C^{*}_{r}(G)^{*}$. On the other hand (cf. \cite{J}), for
each sequence
$\{k_{j}\}$ of natural numbers,
\begin{align*}
\lim_{n}\bigg\|\frac{1}{n}\sum_{j=1}^{n}\a^{k_{j}}(\l_{g})\bigg\|
\leq&
C(1+L(g))^{s}\bigg\|\frac{1}{n}\sum_{j=1}^{n}\d_{\b^{k_{j}}(g)}\bigg\|_{\ell^{2}(G)}\\
\equiv&\frac{C(1+L(g))^{s}}{\sqrt{n}}\,.
\end{align*}
The assertion follows by Theorem \ref{z}
\end{pf}

Finally, we have the case of the automorphism generated by the shift on
the
free group on infinitely many generators.
\begin{cor}
Let $\bbf_{\infty}$ be the free group on infinitely many generators
$\{g_{i}\}_{i\in\bz}$. The automorphism $\a$ induced on
$C^{*}_{r}(\bbf_{\infty})$
by the free shift of the generators is $E$--strictly weak mixing with
$E=\t(\,\cdot\,)\idd$, $\t$ being the canonical trace on
$C^{*}_{r}(\bbf_{\infty})$.\footnote{The $C^{*}$--dynamical system
$\big(C^{*}_{r}(\bbf_{\infty}),\a,\t\big)$ is indeed {\it strictly weak
mixing} in the
language of \cite{MT}.}
\end{cor}
\begin{pf}
By taking into account the Haagerup inequality (cf. \cite{H}, Lemma
1.4), the
situation under consideration is a particular case of
Proposition \ref{ue1}.
\end{pf}

\section{strict weak mixing of the free-shift on the reduced free product
$C^{*}$--algebras}

The present section is devoted to show that the shift on the reduced
amalgamated free product $C^{*}$--algebra
${\displaystyle (A,\phi)={(*_B)}_{i\in I}(A_i,\phi_i)}$
is indeed $\phi$--strictly weak
mixing, the last being stronger than unique ergodicity w.r.t. its
fixed subalgebra, proven in
\cite{AD}. Also this proof relies on the analogue of the Haagerup
inequality proven in \cite{AD} (cf. Proposition 5.1).

We briefly recall some facts on the reduced
amalgamated free product. For an exhaustive treatment of the
subject see \cite{AD, V}.

Let $D$ be a unital $C^{*}$--algebra with identity $\idd$, and
$E^D_B:D\mapsto B$
a conditional expectation onto the unital $C^{*}$--subalgebra $B$ with
the same identity $\idd$. For
each integer $i\in\bz$ consider a copy $(A_i,\phi_i)$ of
$(D,E^D_B)$, together with the reduced amalgamated free product
\begin{equation}
\label{eq:Aphi}
(A,\phi)={(*_B)}_{i\in I}(A_i,\phi_i)\,.
\end{equation}
The $C^{*}$--algebra $A$ naturally acts on a Hilbert right $B$--module $E$
and it is
generated by $\big\{\l^{i}_{a}\,:\, a\in A_{i}\,,i\in\bz\big\}$,
$\l^{i}$ being the embedding of $A_{i}$ in $\cb_{B}(E)$, the space of
all the bounded $B$--linear maps acting on $E$. The
conditional expectation $\phi$ is given by
$$
\phi(a)=\langle\idd a\,,\,\idd\rangle\,,\qquad a\in A\,,
$$
$\langle\,\cdot\,,\,\cdot\,\rangle$
being the $B$--valued inner product of $E$ which is supposed be
linear w.r.t. the first variable.\footnote{Relatively to the Hilbert
$C^{*}$--modules, see e.g. \cite{Pa}.} The free--shift automorphism
$\alpha$ on $A$ is the automorphism of $A$ given by
$\alpha(\lambda_a^i)=\lambda_a^{i+1}$ for all $a\in A$ and
$i\in\bz$.
\begin{thm}
\label{F-mix}
Let $\alpha$ be the free--shift automorphism on the reduced amalgamated
free product $C^*$--algebra $A$ given in \eqref{eq:Aphi}.
Then $\alpha$ is $\phi$--strictly weakly mixing.
\end{thm}
\begin{pf}
It was proven in \cite{AD} that $\a$ is
$\phi$--uniquely ergodic, that is
\begin{equation*}
\lim_{n\to+\infty}\frac1n\sum_{k=0}^{n-1}\alpha^k(a)=\phi(a)
\end{equation*}
for every $a\in A$. By a standard density argument, it is enough to prove
that the sequence $\{\a^{n}(a)\}_{n\geq 1}$ is weakly mixing to
zero whenever $a$ has the form $a=w$ for a word
$w=\l^{m(1)}_{a_1}\l^{m(2)}_{a_2}\cdots\l^{m(p)}_{a_p}$,
with $p\geq 1$, $a_i\in A_{m(i)}^\circ$,
and $m(i)\in\bz$ fulfilling $m(i)\ne m(i+1)$,
$i=1,\dots,p-1$. Here,
$$
A_{i}^\circ:=\big\{a-\phi_{i}(a)\,:\,a\in A_{i}\big\}\,,\qquad
i\in\bz\,.
$$
Let us take any increasing sequence
$\{k_j\}\subset\bn$.
Notice that
$$
\a^{k}(w)=\l^{m(1)+k}_{a_1}\l^{m(2)+k}_{a_2}\cdots\l^{m(p)+k}_{a_p}\,,
$$
that is $\a^{k}(w)$ is a word satisfying the same properties as $w$
with $m'(i)=m(i)+k$.\footnote{This simply means that
$$
\phi\big(w^{*}\a^{k}(w)\big)=\d_{k,0}(2p+1)^{2}\prod_{i=1}^p\|a_i\|^{2}
$$
or equivalently, $\big\{\idd\a^{k}(w)\big\}_{k\in\bz}\subset E$ is an
orthogonal set in the $B$--module $E$.}
Then we can apply the estimation in Proposition 5.1 of \cite{AD} to
the element
$$
f:=\sum_{j=1}^{n}\a^{k_{j}}(w)
$$
obtaining
$$
\bigg\|\frac{1}{n}\sum_{j=1}^{n}\alpha^{k_j}(w)\bigg\|
\leq\frac{2p+1}{n^{1/2}}\prod_{i=1}^p\|a_i\|\,.
$$
Now, by applying Theorem \ref{z}, we obtain the assertion.
\end{pf}

\section*{acknowledgments} The second named--author (F. M.) thanks FCT
(Portugal) grant SFRH/ BPD/17419/2004, and Prof. L. Accardi for
kind hospitality  from 18--22 June 2006 at ``Universit\`{a} di Roma
Tor Vergata''.

\end{document}